\documentclass[12pt,a4paper]{amsart}
\usepackage{amsmath,amsfonts}
\usepackage{graphicx}
\usepackage{epsfig}
\usepackage{amssymb}
\usepackage{hyperref}

\newtheorem{theorem}{Theorem}[section]

\newtheorem{defprop}[theorem]{Definition and Proposition}

\newtheorem*{decotheorem}{Decoration Theorem}
\newtheorem*{maintheorem}{Main Theorem}

\theoremstyle{definition}

\newtheorem*{remark*}{Remark}

\newcommand{\CC}{\mathbb{C}}

\newcommand{\LL}{\mathcal L}
\newcommand{\DD}{\mathcal D}
\newcommand{\MM}{\mathcal M}

\newcommand{\HH}{\mathcal H}
\newcommand{\cHH}{\widehat{\HH}}
\newcommand{\disk}{\mathbb D}

\newcommand{\eps}{\varepsilon}
\newcommand{\sm}{\setminus}
\newcommand{\ovl}{\overline}

\begin{document}

\title[The Homeomorphism Theorem ]{Homeomorphisms between limbs of the Mandelbrot set}
\author{Dzmitry Dudko}
\address{Research I, Jacobs University, Postfach 750 561, D-28725 Bremen, Germany; and: G.-A.-Universit\"at zu G\"ottingen, Bunsenstrasse 3--5, D-37073 G\"ottingen, Germany}
\email{d.dudko@jacobs-university.de}
\author{Dierk Schleicher}
\address{Research I, Jacobs University, Postfach 750 561, D-28725 Bremen, Germany}
\email{dierk@jacobs-university.de}
\date{\today}

\begin{abstract}
We prove that for every hyperbolic component of the Mandelbrot set, any two limbs with equal denominators are homeomorphic so that the homeomorphism preserves periods of hyperbolic components. This settles a conjecture on the Mandelbrot set that goes back to 1994.
\end{abstract}

\maketitle

\section{Introduction}

The Mandelbrot set $\MM$ is a set with a very rich combinatorial, topological, and geometric structure. It is often called ``self-similar'' because there are countably many dynamically defined homeomorphisms from $\MM$ into itself, and the set of such homeomorphisms forms a semigroup. Moreover, there are many dynamically defined homeomorphisms from certain dynamically defined subsets of $\MM$ to other subsets of $\MM$. Perhaps the first such result was a homeomorphism from the $1/2$-limb of $\MM$ to a subset of the $1/3$-limb of $\MM$ constructed by Branner and Douady \cite{BD}; this class of homeomorphisms was later extended by Riedl \cite{Ri}.

In \cite{BF1}, it was shown, using homeomorphisms to parameter spaces of certain higher degree polynomials, that any two limbs $\LL_{p/q}$ and
$\LL_{p'/q}$ (with equal denominators) were homeomorphic. These homeomorphisms preserve the embedding into the plane so that they even extend to neighborhoods of these limbs within $\CC$, preserving the orientation \cite{BF2}. All these homeomorphisms are constructed by quasiconformal surgery, and they all change the dynamics of the associated polynomials so that in general, periods of hyperbolic components are changed.

At about the same time, it was observed \cite{LS} that there is a combinatorially defined bijection between the limbs $\LL_{p/q}$ and
$\LL_{p'/q}$ that preserves periods of hyperbolic components, and it was conjectured that this would yield a homeomorphism between these limbs that preserved periods of hyperbolic components. An early attempt to prove this conjecture by quasiconformal surgery resulted in another proof of the theorem from \cite{BF1} that stayed within the quadratic family.

A proof of this conjecture is the main result of the present paper; it can be stated as follows.
\begin{maintheorem}
For any hyperbolic component of $\MM$, let $\LL_{p/q}$ and $\LL_{p'/q}$ be two limbs with equal denominators. Then there exists a homeomorphism between them that preserves periods of hyperbolic components.
\end{maintheorem}

Since our homeomorphism preserves periods of hyperbolic components, it can not
 extend to neighborhoods of the limbs.

\begin{figure}[hb]
\includegraphics[width=6cm]{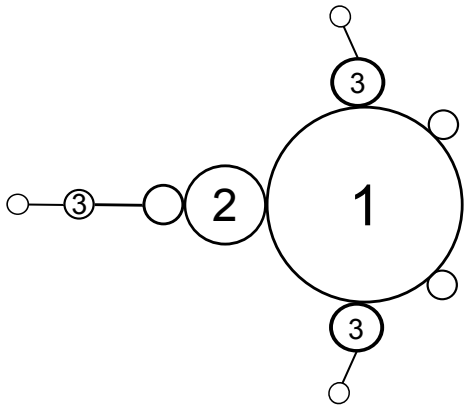}
  \caption{Combinatorics of hyperbolic components of $\MM$ up to period 4.
  }
  \label{figure:CombMand}
\end{figure}

For a fixed $n\ge1$ consider the arrangement $\MM_n$ of all
hyperbolic components with periods up to $n$ (see Figure
\ref{figure:CombMand} for an example). There is a combinatorial
model $\MM_{comb}$ of the Mandelbrot set that can be described as
a limit of $\MM_n$ in a certain sense \cite{Do}. Furthermore,
there is a canonical continuous projection $\pi:\MM\rightarrow
\MM_{comb}$, and any fiber $\pi^{-1}(c)$ is compact, connected,
and full (a bounded set $X\subset\CC$ is called \emph{full} if its
complement has no bounded components). The famous ``MLC
conjecture'' (``the Mandelbrot set is locally connected'') can be
stated as saying that  $\pi$ is a homeomorphism.

For any $p/q$ and $p'/q$ there is a canonical homeomorphism $f'$
between $\pi(L_{p/q})$ and $\pi(L_{p'/q})$ preserving periods of hyperbolic components. Our strategy is to
show that $f'$ can be lifted up to the level of the Mandelbrot set;
namely we have the following commutative diagram:
\begin{equation}
\label{eq:diagram0}
\begin{array}[c]{ccc}
L_{p/q}&\stackrel{f}{\longrightarrow} &L_{p'/q}\\
\downarrow\scriptstyle{\pi}&&\downarrow\scriptstyle{\pi}\\
\pi(L_{p/q})&\stackrel{f'}{\longrightarrow}&\pi(L_{p'/q}).
\end{array}
\end{equation}
We will show that this technique can be applied to any continuous map that
``respects'' small copies of the Mandelbrot set.

This result fits into the vision of Douady expressed by the statement that ``combinatorics implies topology'': many results about the Mandelbrot sets are discovered and described in terms of combinatorics, and these combinatorial results lead the way for topological statements. In our case, the combinatorial result remained a topological conjecture since about 1994. The key progress that was required was the Decoration Theorem (see below).

\subsection*{Outline of the paper.} In Section \ref{sect:Mand} we
recall the notion of hyperbolic components, small copies of the
Mandelbrot set, and combinatorial classes. The combinatorial model
$\MM_{comb}$ is defined as the quotient of $\MM$.

Section \ref{sect:internal_addresses} contains the definition and
main properties of internal and angled internal addresses. They
are coordinates for combinatorial classes.

In Section \ref{sec:proof} we will construct Diagram
\ref{eq:diagram0}. The homeomorphism $f':\pi(L_{p/q})\rightarrow
\pi(L_{p'/q})$ exists by fundamental properties of angled internal
addresses. As $f'$ coincides with the canonical homeomorphism on
every small copy of the Mandelbrot set there exists a bijection
$f$ that makes Diagram \ref{eq:diagram0} commute. The continuity
of $f$ follows from Yoccoz's results, the existence of the canonical
isomorphism of all copies of the Mandelbrot set, and the
Decoration theorem.

In Section \ref{sec:generalization} we will formulate a general
statement that allows to lift a continuous map from the level of
the combinatorial model of the Mandelbrot set to the actual
Mandelbrot set.

\medskip





\section{The Mandelbrot set}
\label{sect:Mand}

The \textit{Mandelbrot set} $\MM$ is defined as the set of
quadratic polynomials
\begin{equation}
\label{eq:Qadr} p_c\colon z\mapsto z^2+c
\end{equation}
 with connected Julia sets. It is a compact,
connected, and full set.

As the Mandelbrot set is the parameter space of quadratic
polynomials there is an additional structure (combinatorics) of
$\MM$ on top of the topology. For instance, $\MM$ contains
hyperbolic components and small copies of the Mandelbrot set. Both
types of subsets have dynamical meaning.

\subsubsection*{Hyperbolic components} A \emph{hyperbolic component of $\MM$} is a connected component of the set of parameters $c\in\MM$ so that $p_c$ has an attracting orbit.

Assume that $p_c(z)=z^2+c:\CC\rightarrow\CC$
has an non-repelling periodic cycle; this means there is a periodic orbit $z_c$ of $p_c$ with multiplier of absolute value at most $1$. The periodic orbit $z_c$ is necessarily unique; let $\lambda(z_c)$ be its multiplier.
Then it is known that:
\begin{itemize}
\item
there is a hyperbolic component $\HH$ such that
$c\in\overline{\HH}\subset\MM$;
\item
the attracting orbit $z_c$ has constant period throughout $\HH$;
\item
within $\HH$ the cycle $z_c$ moves holomorphically;
\item the multiplier map $\lambda\colon\HH\to\disk$
 is a conformal isomorphism; this map extends to a homeomorphism from $\ovl \HH$ to $\ovl\disk$.
\end{itemize}

By definition, the period of $\HH$ is the period of $z_c$. For
every fixed $n\ge 1$ there are finitely many hyperbolic components
with period $n$. The arrangement of hyperbolic
components up to period $n$ gives an approximation (topological
and combinatorial) to the Mandelbrot set (see Figure
\ref{figure:CombMand}). The unique hyperbolic component of period
$1$ is called the \textit{main hyperbolic component of the
Mandelbrot set}.

Let $\HH_1$ and $\HH_2$ be hyperbolic components with periodic
orbits $z_c$ and $z'_c$ and periods $n_1$ and $n_2$,
respectively. If the closures of $\HH_1$ and $\HH_2$ intersect,
then the intersection is one point. Let us assume that $\HH_1$ and
$\HH_2$ intersect, $n_1\le n_2$, and the point $c'$ is the
intersection.

It is known that:

\begin{itemize}

\item $n_2/n_1=q$ is an integer greater than $1$;

\item at $c'$ the cycles $z_c$ and $z'_c$ collide: every point
from the cycle $z_c$ collides with $q$ points from the cycle
$z'_c$;

\item the multiplier of $z_c$ at $c'$ is $\exp(2\pi i p/q )$, where $p$ is
an integer coprime (and less than) to $q$; in particular, the
unique non-repelling orbit of $p_{c'}$ is parabolic.
\end{itemize}
The closure of the connected component of $\MM\setminus\HH_1$
containing $\HH_2$ is called the \emph{$p/q$-limb}
$\LL_{p/q}(\HH_1)$ \textit{of} $\HH_1$. It is known that
$p/q$-limbs exist for all coprime $p/q$ with $q\ge 2$,
 and (the closure of) every component of $\MM\sm\HH_1$ that does not contain the point $c=0$ is a limb of $\HH_1$.
If $\HH_1$ is the main hyperbolic component, then
$\LL_{p/q}(\HH_1)$ is called the (primary) \textit{$p/q$-limb}
$\LL_{p/q}$.

We define \textit{the combinatorial class} $\cHH$ of a
hyperbolic component $\HH$ as $\{c\in\overline{\HH}\ |\
\lambda(z_c)\not= p/q,\ q>1\}$, where $\lambda(z_c)$ is the
multiplier of the (unique) non-repelling periodic cycle $z_c$ of a
polynomial $p_c$; equivalently:
\begin{equation}
\label{eq:CombClos}
\cHH=\overline{\HH}\setminus\bigcup_{q=2}^{\infty}\bigcup_{p}\LL_{p/q}(\HH).
\end{equation}

 \subsubsection*{Small copies of the Mandelbrot set} The
 Mandelbrot set is a self-similar set in the following sense: there are countably
 many copies of the Mandelbrot set in $\MM$; every such copy is canonically
 homeomorphic to $\MM$, where the homeomorphism is given by the straightening theorem \cite{DH}. In particular small copies are compact, connected, full sets and they preserve
 hyperbolic components: if a small copy intersects a hyperbolic component, then it contains it.
Polynomials within small copies of $\MM$ are called \emph{renormalizable}; polynomials within infinitely many nested small copies of $\MM$ are called \emph{infinitely renormalizable}.

 Small copies are in one-to-one
 correspondence with hyperbolic components: for every hyperbolic component $\HH$, there is a unique small copy $\MM_\HH\supset\HH$ so that the canonical homeomorphism of $\MM_\HH$ sends $\HH$ to the main hyperbolic component of $\MM$, and every small copy of $\MM$ is of this type for a unique component $\HH$.
For a small copy $\MM_\HH$, let $per(\MM_\HH)$ be the period of $\HH$. Then the canonical homeomorphism from $\MM_\HH$ to $\MM$ divides all periods of hyperbolic components by $per(\MM_\HH)$.

If $\MM_\HH$ is a small copy of $\MM$ within $\MM$, then $\MM\sm\MM_\HH$ consists of countably many components, called ``decorations'' of $\MM_\HH$. If $\DD$ is the closure of any such decoration, then $\DD\cap\MM_\HH$ is a single Misiurewicz point (i.e., a parameter for which the critical orbit is strictly preperiodic). The following  theorem was recently proved \cite{D} (see \cite{PR} for a different proof); it will be the fundamental motor for our theorem.

\begin{decotheorem}
 For any $\varepsilon >
0$, there are at most finitely many connected components of $\MM
\backslash \MM_{s}$ with diameter at least $\eps$.
\end{decotheorem}

\subsubsection*{Yoccoz polynomials} A \emph{Yoccoz polynomial} is a quadratic polynomial in $\MM$ for which all periodic orbits are repelling, and it is not infinitely renormalizable
(equivalently, it does not belong to any hyperbolic component and
is not within infinitely many small copies of $\MM$). It is known
that $\MM$ is locally connected at Yoccoz polynomials, and
stronger yet, that the corresponding fibers of $\MM$ are trivial:
this was shown in detail in \cite{Hu} for non-renormalizable
parameters, but results of this kind are automatically preserved
by finite renormalizations \cite{Sch2}.

\subsection*{Combinatorial classes}
A combinatorial class is an equivalence class of parameters with
the same rational lamination. A \textit{combinatorial class} is either:

\begin{itemize}
\item
a hyperbolic combinatorial class (as defined above); or
\item
the intersection of an infinite nested sequence of small copies of
$\MM$;

\item
a single point that does not belong to any combinatorial class of the
first two types.
\end{itemize}
Any
non-hyperbolic combinatorial class is always a compact, connected,
and full set. A combinatorial class of the last type is exactly a
Yoccoz parameter.

 There are two famous conjectures: ``The
Mandelbrot set is locally connected'' (MLC) and ``hyperbolic
dynamics is dense in the space of quadratic polynomials'' (the
Fatou conjecture for quadratic polynomials). These are equivalent
to the statements  ``every non-hyperbolic combinatorial class is a
point'' and ``every non-hyperbolic combinatorial class has no
interior point'' respectively \cite{Do}.

\subsection*{The combinatorial model of the Mandelbrot set.}
Let us say that two points $c_1$ and $c_2$ are \emph{combinatorially
equivalent} if $c_1$ and $c_2$ are in the same non-hyperbolic
combinatorial class. The \textit{combinatorial model} $\MM_{comb}$
of the Mandelbrot set is the quotient of $\MM$ by the above
equivalence relation. The associated projection $\pi:\MM
\rightarrow\MM_{comb}$ is called canonical.

It is known that $\MM_{comb}$ is a connected, locally connected,
compact, full set; $\pi$ is a continuous surjection and $\pi$ is a
homeomorphism (i.e., injective) if and only if MLC is valid \cite{Do}.

Hyperbolic components and small copies for $\MM_{comb}$ are
defined using the projection $\pi$.

Yoccoz' theorem can be expressed as follows: if $c_n$ is a sequence of parameters in $\MM$ so that $\pi(c_n)\to\pi(c)$ for some Yoccoz parameter $c\in\MM$, then $c_n\to c$; and this is the statement we need. (Note that this does \emph{not} follow from the fact that $\MM$ is locally connected at $c$; the stronger property is required that the fiber of $\MM$ at $c$ is trivial; and Yoccoz indeed proves that; see \cite{Sch2}.)

\section{Internal addresses of the Mandelbrot set}
\label{sect:internal_addresses}

In this section we recall the definition and main properties of
internal and angled internal addresses. The main reference is
\cite{Sch1}. The motivation of an internal address is
  to approximate any combinatorial class by a canonical sequence of (simpler) hyperbolic
 classes. Internal addresses (and angled internal addresses) are defined for
 combinatorial classes, hence there is no difference in
 the definitions for $\MM$ and $\MM_{comb}$.

Consider a combinatorial class $C$ and a hyperbolic component
$\HH$. Assume that either $\cHH=C$ or $C$ is not in the connected
component of $\MM\setminus C$ containing $0$. Then we say that
$\HH$ is \textit{closer to} $0$ than $C$ and write $\HH\le C$. We
also write $C<\HH$ if $C\le\HH$ and $\cHH\not=C$.

 For any $C$ inductively define the (finite
or infinite) sequence
\begin{equation}
\label{eq:IntAddr}  \HH_0< \HH_1<\dots<\HH_n< \dots
\end{equation}
 such that $\HH_0$ is the main hyperbolic component of $\MM$
and $\HH_n$ is of the smallest period satisfying $\HH_{n-1}< \HH_n
\le C$.

\begin{defprop}
For any $C$ the sequence in (\ref{eq:IntAddr}) is unique. Define $S_n$ to
be the period of $\HH_n$ and let $p_n/q_n$ be the fraction so that $\HH_{n+1}\subset\LL_{p_n/q_n}(\HH_n)$.

The sequence
\begin{equation}1=S_0 \rightarrow S_1\rightarrow S_2
\rightarrow \dots
\end{equation} is called the \textbf{internal}
address of $C$.

The sequence
\begin{equation} (S_0)_{p_0/q_0} \rightarrow (S_1)_{p_1/q_1}\rightarrow
(S_2)_{p_2/q_2} \rightarrow \dots
\end{equation} is called the
\textbf{angled internal} address of $C$.
\end{defprop}

   It is known \cite[Theorem 1.10]{Sch1} that an angled internal address uniquely
   describes a combinatorial class, where finite addresses correspond to hyperbolic classes.
    On the other hand the
   internal address describes a combinatorial class up to the ``symmetry''.
   For example, hyperbolic polynomials have the same internal addresses if and only if the dynamics of the polynomials
    on the Julia sets are topologically conjugate; this topological conjugation extends to a neighborhood of the Julia set, preserving the orientation, if and only if the two polynomials have the same angled internal address.

Internal addresses are strictly increasing (finite or infinite) sequences of integers starting with $1$. Not every such sequence occurs for a combinatorial class of the Mandelbrot set: those that occur are called ``complex admissible''. An explicit characterization of the complex admissible sequences was given in \cite{BS}.

    The   following theorem shows the
   ``valency of the symmetry'':

\begin{theorem}[{\cite[Theorem 2.3]{Sch1}}]
\label{th:AnglIntAddr} If an angled internal address describes a
combinatorial class in the Mandelbrot set, then the numerators
$p_k$ can be changed arbitrarily (coprime to $q_k$) and the
modified angled internal address still describes a combinatorial
class in the Mandelbrot set.
\end{theorem}

In other words, complex admissibility is a property of internal addresses, not of angled internal addresses. In fact, any internal address uniquely determines the denominators $q_k$ of any associated angled internal address, while Theorem~\ref{th:AnglIntAddr} says that the numerators $p_k$ are completely arbitrary (coprime to $q_k$).

Consider a small copy $\MM'$ of the Mandelbrot set, and assume
$\HH'$ is the main hyperbolic component of $\MM'$. Then $\HH'$ has
a finite angled internal address:
\begin{equation} (S'_0)_{p'_0/q'_0} \rightarrow (S'_1)_{p'_1/q'_1}\rightarrow
(S'_2)_{p'_2/q'_2} \rightarrow \dots \rightarrow S'_n.
\end{equation}

\begin{theorem}[{\cite[Proposition 2.7]{Sch1}}]
\label{th:SmallCop} A combinatorial class $C$ belongs to the small
copy $\MM'$ if and only if the angled internal address of $C$ is
\begin{equation} (S'_0)_{p'_0/q'_0} \rightarrow
  \dots  (S'_{n-1})_{p'_{n-1}/q'_{n-1}}\rightarrow (S'_n)_{p_n/q_n}\rightarrow (S_{n+1})_{p_{n+1}/q_{n+1}}\rightarrow\dots
\end{equation}
and $S'_n|S_{n+k}$ for $k\ge 1$. The canonical homeomorphism
between $\MM'$ and $\MM$ sends $C$ to the combinatorial class with
the internal address
\begin{equation}  (1)_{p_n/q_n}\rightarrow (S_{n+1}/S'_n)_{p_{n+1}/q_{n+1}}\rightarrow
 (S_{n+2}/S'_n)_{p_{n+1}/q_{n+1}}\rightarrow
\dots.
\end{equation}
\end{theorem}

We need one more result from \cite{Sch1}.
\begin{theorem}[{\cite[Proposition 2.6]{Sch1}}]
\label{th:BranchIntAddr}
Consider a hyperbolic component $\HH$ and a combinatorial class $C$ in the $p/q$-limb of $\HH$. If $q\ge 3$, then $\HH$ occurs in the internal address of $C$; more precisely, the internal address of $\HH$ is a finite initial sequence of the internal address of $C$.
\end{theorem}

This result can be expressed as follows: for a given combinatorial class $C$, there are usually many hyperbolic components $\HH<C$, and most of them are not associated to the internal address of $C$. For those that are not, $C$ is in the $1/2$-limb of $\HH$.

\section{Proof of the homeomorphism theorem}
\label{sec:proof}

In this, section, we prove the main theorem in an apparently stronger form: consider hyperbolic components $\HH_1$ and $\HH_2$ with identical
internal addresses. Then we will construct a
homeomorphism between the limbs $\LL_{p/q}(\HH_1)$ and
$\LL_{p'/q}(\HH_2)$, where $q\ge 3$. The original statement of the Main Theorem describes the case $\HH_1=\HH_2$.

\begin{remark*}
This more general version can easily be deduced from the statement
of the Main Theorem, because there is a unique hyperbolic $\HH'$
at which the angled internal addresses of $\HH_1$ and $\HH_2$
branch off, in the sense that $\HH_1$ and $\HH_2$ are contained in
two different limbs at angles $p/q$ and $p'/q$ of $\HH'$, with
$q\ge 3$; a possibly repeated application of the Main Theorem will
then yield the statement we are proving here, and it shows that
the statement remains true for $q=2$, i.e., for the limbs
$\LL_{1/2}(\HH_1)$ and $\LL_{1/2}(\HH_2)$. The reason why we are
giving an apparently more general proof is because the proof
really is the same, and this illustrates the general nature of the
argument.
\end{remark*}

Let the angled internal addresses of $\HH_1$ and $\HH_2$ be
\begin{equation}
 (S_0)_{p_0/q_0} \rightarrow\dots\rightarrow
(S_{n-1})_{p_{n-1}/q_{n-1}}\rightarrow (S_n) \;,
\end{equation}
\begin{equation}
(S_0)_{p'_0/q_0} \rightarrow\dots\rightarrow
(S_{n-1})_{p'_{n-1}/q_{n-1}}\rightarrow (S_n)
\end{equation}
respectively. By Theorem~\ref{th:BranchIntAddr}, the limbs $\LL_{p/q}(\HH_1)$ and
$\LL_{p'/q}(\HH_2)$ consist exactly of all combinatorial classes that have internal addresses starting with
\begin{equation}
\label{eq:pr1} (S_0)_{p_0/q_0} \rightarrow\dots\rightarrow
(S_{n-1})_{p_{n-1}/q_{n-1}}\rightarrow (S_n)_{p/q}\rightarrow\;,
\end{equation}
\begin{equation}
\label{eq:pr2} (S_0)_{p'_0/q_0} \rightarrow\dots\rightarrow
(S_{n-1})_{p'_{n-1}/q_{n-1}}\rightarrow (S_n)_{p'/q}\rightarrow
\end{equation}
respectively.

Define a map $f'\colon\pi(\LL_{p/q}(\HH_1))\rightarrow
\pi(\LL_{p'/q}(\HH_2))$ so that it changes the initial segment (\ref{eq:pr1}) of the angled internal address into the segment (\ref{eq:pr2}), i.e., it changes the angles in the angled internal address from the limb $\LL_{p/q}(\HH_1)$ into the limb $\LL_{p'/q}(\HH_2)$; within hyperbolic components, this map should fix multipliers. It follows from the construction and Theorem~\ref{th:AnglIntAddr} that the map $f'$ is a well defined homeomorphism and it preserves internal addresses. We will show that there exists a
canonical homeomorphism $f$ such that:
\begin{equation}
\label{eq:diagram}
\begin{array}[c]{ccc}
\LL_{p/q}(\HH_1)&\stackrel{f}{\longrightarrow} &\LL_{p'/q}(\HH_2)\\
\downarrow\scriptstyle{\pi}&&\downarrow\scriptstyle{\pi}\\
\pi(\LL_{p/q}(\HH_1))&\stackrel{f'}{\longrightarrow}&\pi(\LL_{p'/q}(\HH_2)).
\end{array}
\end{equation}
By canonical we mean that $f$ coincides with the natural
homeomorphism between small copies of the Mandelbrot set. With
this requirement there is a unique bijection $f$ that makes the
above diagram commute. Indeed, if a non-hyperbolic combinatorial
class $C$ does not belong to any small copy of $\MM$, then $C$ is
a point and $f(C)$ is uniquely defined. If $C$ belongs to a small
copy $\MM_s\subset\LL_{p/q}(\HH_1)$ of $\MM$, then $f$ on $\MM_s$
is uniquely defined by the requirement that $f$ be canonical (note
that $f'$ is canonical by Theorem \ref{th:SmallCop}).

The main issue is to prove that $f$ is continuous. Let us prove
that if $c_n$ tends to $c_\infty$, then  $f(c_n)$ tends to $f(c_\infty)$.
This will imply that $f$ is a homeomorphism (as $f$ is a continuous bijection between compact Hausdorff spaces). It suffices to consider the following three cases.

\subsection*{Case 1} Assume $c_\infty$ belongs to at most finitely many small
copies of the Mandelbrot set; then the same is true for $f(c_\infty)$. By construction, $f'(\pi(c_n))=\pi(f(c_n))$ tends to $f'(\pi(c_\infty))=\pi(f(c_\infty))$ (using commutativity of Diagram
(\ref{eq:diagram})). By Yoccoz' theorem, it follows that $f(c_n)$ tends to $f(c_\infty)$.

\subsection*{Case 2} Assume $c_\infty$ and all $c_n$ belong to a single small copy $\MM_s$, where $\MM_s\subset\LL_{p/q}(\HH_1)$.
Then the statement follows from the definition of $f$ because $f$
coincides  with the canonical homeomorphism from $\MM_s$ to
$f(\MM_s)$.

\subsection*{Case 3} Assume $c_\infty$ belongs to infinitely many copies
of the Mandelbrot set (i.e., $c_\infty$ is infinitely
renormalizable), $\MM_s\subset\LL_{p/q}(\HH_1)$ is a small copy
containing $c_\infty$, and
 $c_n$ does not belong to $\MM_s$ for any $n$.

 Let $\DD_n$ be the closure of the connected component of $\MM
\backslash \MM_{s}$ containing $c_n$, and let $a_n$ be the
intersection of $\MM_{s}$ and $\DD_n$. Then $a_n$ is a single Misiurewicz point and hence belongs to at most finitely many copies of $\MM$.

 As $c_\infty$ belongs to infinitely many copies of $\MM$, it follows that
 $c_\infty\neq a_k$ for all $k$. Therefore only finitely many $c_n$ are in $\DD_k$ for each fixed $k$. Hence by the Decoration Theorem the distance between $c_k$ and $a_k$ tends to
 $0$, and so
 the sequence $a_k$ tends to $c_\infty$.
  By Case $2$ we obtain that $f(a_k)$ tends to
 $f(c_\infty)$.

 By a similar reason the distance between $f(a_k)$ and $f(c_k)$
 tends to $0$ ($f(a_k)$ is the intersection of $f(\MM_s)$ with the closure of the
  connected component of $\MM\setminus \MM_s$ containing
  $f(c_k)$). We conclude that $f(c_k)$ tends to $f(c_\infty)$. This concludes the proof of the Main Theorem.

\section{Generalization}
\label{sec:generalization}

  We say that a set $\LL\subset\MM$ is \textit{combinatorially saturated} if
  $\pi^{-1}(\pi(\LL))=\LL$.
  Consider two closed combinatorially saturated sets $\LL_1$ and $\LL_2$ and
  assume that there exists a continuous map $f':\pi(\LL_1)\rightarrow
  \pi(\LL_2)$.

   We say that $f'$ is \textit{canonical} (with respect to small copies of $\MM$) if:
 \begin{itemize}
   \item for every infinitely
renormalizable $c\in\pi(\LL_1)$ there exists a copy
$\pi(\MM_s)\subset\pi(\LL_1)$ containing $c$ such that $f'$ restricted to
$\pi(\MM_s)$ is the canonical homeomorphism on $\pi(\MM_s)$;

   \item for every infinitely
renormalizable $c\in\pi(\LL_2)$ there exists a copy
$\pi(\MM_s)\subset\pi(\LL_2)$ containing $c$ such that $f'$
restricted to any connected component of $f'^{-1}(\pi(\MM_s))$ is
the canonical homeomorphism of a small copy of the Mandelbrot set.
\end{itemize}
In particular, by a standard compactness argument
$f'^{-1}(\pi(\MM_s))$ (in the second condition) consists of
finitely many small copies.

\begin{theorem}
Let $\LL_1,\LL_2\subset\MM$ be two closed connected combinatorial
sets.
 For every continuous map $f':\pi(\LL_1)\rightarrow
  \pi(\LL_2)$ that is canonical with respect to small copies of the Mandelbrot set
  there exists a continuous map $f:\LL_1\rightarrow
  \LL_2$ such that the following diagram commutes:
\begin{equation}
\begin{array}[c]{ccc}
\LL_{1}&\stackrel{f}{\longrightarrow} &\LL_{2}\\
\downarrow\scriptstyle{\pi}&&\downarrow\scriptstyle{\pi}\\
\pi(\LL_{1})&\stackrel{f'}{\longrightarrow}&\pi(\LL_{2}).
\end{array}
\end{equation}
\end{theorem}
The proof is quite similar to previous one
and is left to the reader.

%
%
%

\end{document}